\documentclass[12pt]{amsart}
\usepackage{a4}
\usepackage{amssymb,amsmath}
\newtheorem{thm}{Theorem}

\newtheorem{cor}[thm]{Corollary}
\newtheorem{lem}[thm]{Lemma}
\newtheorem{prop}[thm]{Proposition}
\theoremstyle{definition}

\theoremstyle{remark}
\newtheorem{rem}{Remark}

\begin{document}
\title{On the Vector valued Fourier Transform And Compatibility of Operators* }
\author{In Sook Park}
\address{Division of Applied Mathematics, Korea Advanced
Institute of Science and Technology, 373-1 Kusong-dong, Yusong-gu,
Taejon 305-701, Republic of Korea}

\email{ispark@amath.kaist.ac.kr}

\thanks{*This work is supported by BK21 project.}
\begin{abstract}
Let $\mathbb{G}$ be a locally compact abelian group and let
$1<p\leq 2$. $\mathbb{G}^{'}$ is the dual group of $\mathbb{G}$,
and $p^{'}$ the conjugate exponent of $p$. An operator $T$ between
Banach spaces $X$ and $Y$ is said to be compatible with the
Fourier transform $F^{\mathbb{G}}$ \,if\, $F^{\mathbb{G}}\otimes
T: L_p(\mathbb{G})\otimes X\rightarrow
L_{p^{'}}(\mathbb{G}^{'})\otimes Y $ admits a continuous extension
$[F^{\mathbb{G}},T]:[L_p(\mathbb{G}),X]\rightarrow
[L_{p^{'}}(\mathbb{G}^{'}),Y]$. $\mathcal{FT}_p^{\mathbb{G}}$
denotes the set of such $T$'s. We show that
$\mathcal{FT}_p^{\mathbb{R}\times\mathbb{G}}=\mathcal{FT}_p^{\mathbb{Z}\times\mathbb{G}}
=\mathcal{FT}_p^{\mathbb{Z}^n \times\mathbb{G}}$ for any
$\mathbb{G}$ and positive integer $n$. And if the factor group of
$\mathbb{G}$ with respect to its component of the identity element
is a direct sum of a torsion free group and a finite group with
discrete topology then
$\mathcal{FT}_p^{\mathbb{G}}=\mathcal{FT}_p^{\mathbb{Z}}$ .

\end{abstract}
\maketitle
\section{Introduction}
A locally compact abelian group means a topological abelian group
with its Haar measure whose topology is locally compact Hausdorff.
The real line $\mathbb{R}$, the discrete group of integers
$\mathbb{Z}$ and the circle group $\mathbb{T}$ are important
examples. Further information can be found in \cite{Fol},
\cite{Hw} and \cite{Rud}. Let $\mathbb{G}$ be a locally compact
abelian group and $\mathbb{G}^{'}$ its dual(charactor) group, the
Haar measure of $\mathbb{G}^{'}$ is determined so that Parseval's
identity is established with constant 1. For $1\leq r < \infty$,
we denote by $[L_r(\mathbb{G},\mu_{\mathbb{G}}),X]$ the Banach
space of all measurable functions
$\mathbf{f}:\mathbb{G}\rightarrow X $ for which
$\|\mathbf{f}|L_r(\mathbb{G})\|:=\Big{(}\int_{\mathbb{G}}\|\mathbf{f}(s)\|^r\,d\mu_{\mathbb{G}}(s)\Big{)}^{1/r}$
is finite. Let $\mathbb{G}$ be a fixed infinite locally compact
abelian group and $1 < p \leq 2$. $F^{\mathbb{G}}$ denotes the
Fourier transform from $L_p(\mathbb{G})$ into
$L_{p^{'}}(\mathbb{G}^{'})$. For a bounded linear operator $T$
between Banach spaces $X$ and $Y$, $$F^{\mathbb{G}}\otimes
T:\sum_{k=1}^n f_k\otimes\mathbf{x}_k\mapsto
\sum_{k=1}^nF^{\mathbb{G}}f_k\otimes T\mathbf{x}_k$$ yields a
well-defined map from $L_p(\mathbb{G})\otimes X$ into
$L_{p^{'}}(\mathbb{G}^{'})\otimes Y $. $T$ is said to be
compatible with $F^{\mathbb{G}}$ or have $\mathbb{G}$-Fourier type
$p$, if the operator $F^{\mathbb{G}}\otimes
T:L_p(\mathbb{G})\otimes X\rightarrow
L_{p^{'}}(\mathbb{G}^{'})\otimes Y $ admits a continuous extension
$$[F^{\mathbb{G}},T]:[L_p(\mathbb{G}),X]\rightarrow
[L_{p^{'}}(\mathbb{G}^{'}),Y].$$ For such $T$, we let
$$\|\,T\,|\mathcal{FT}_p^{\mathbb{G}}\|:=\|[F^{\mathbb{G}},T]:[L_p(\mathbb{G}),X]\rightarrow[L_{p^{'}}(\mathbb{G}^{'}),Y]\|.$$
The norm defined above is invariant under changing the Haar
measure of $\mathbb{G}$. The class of these operators is a Banach
ideal, denoted by $\mathcal{FT}_p^{\mathbb{G}}$. The definition
and notation follow those of \cite{Pie}.

It is known that
$\mathcal{FT}_p^{\mathbb{R}}=\mathcal{FT}_p^{\mathbb{Z}}=\mathcal{FT}_p^{\mathbb{T}}$
\cite{Pie}, but the problem whether the operator ideal
$\mathcal{FT}_p^{\mathbb{G}}$ depends on $\mathbb{G}$ or not is
unsolved. There are several results about
$\mathrm{FT}_p^{\mathbb{G}}$, the class of Banach
 spaces whose identity operators are compatible with $F^{\mathbb{G}}$.
 And these are immediately extended to the case of $\mathcal{FT}_p^{\mathbb{G}}$ by replacing the identity operator on
 a Banach space $X$ with $T$. Peetre \cite{Pee} who introduced the concept of Banach space of Fourier type $p$ proved that
 $X$ belongs to $\mathrm{FT}_p^{\mathbb{R}}$ if and only if the dual space $X^{'}$ belongs to
 $\mathrm{FT}_p^{\mathbb{R}}$. In fact, $T$ belongs to
 $\mathcal{FT}_p^{\mathbb{G}}$ if and only if the dual operator $T^{'}$ belongs to
 $\mathcal{FT}_p^{\mathbb{G}^{'}}$ i.e.
 $\|\,T\,|\mathcal{FT}_p^{\mathbb{G}}\|=\|\,T^{'}\,|\mathcal{FT}_p^{\mathbb{G}^{'}}\|$
 for any locally compact abelian group $\mathbb{G}$.
 Bourgain \cite{Br} showed that $\mathrm{FT}_p^{\mathbb{T}}
 \subset \mathrm{FT}_p^{\mathbb{R}}$ and K\"{o}nig \cite{Ko} modified Kwapien's argument \cite{Kw} to show that
 $\mathrm{FT}_p^{\mathbb{R}} = \mathrm{FT}_p^{\mathbb{T}}$ and extended this to that
 $\mathrm{FT}_p^{\mathbb{G}} = \mathrm{FT}_p^{\mathbb{T}}$ if $\mathbb{G}$ is one of
 $\mathbb{R}^m$ and $\mathbb{T}^m$, where $m$ is a positive
 integer. Garcia-Cuerva, Kazarian and Torrea \cite{Gar} and
 Andersson \cite{An} showed independently that $\mathrm{FT}_p^{\mathbb{G}} =
 \mathrm{FT}_p^{\mathbb{Z}}$ whenever $\mathbb{G}$ is one of $\mathbb{T}^m,
 \ \mathbb{T}^{\infty},\ \mathbb{R}^m,\ \mathbb{Z}^m$ and $\mathbb{Z}^{\infty}$.
 Andersson \cite{An} also proved that
 $\|I_X|\mathcal{FT}_p^{\mathbb{H}}\|\leq\|I_X|\mathcal{FT}_p^{\mathbb{G}}\|$
  when $\mathbb{H}$ is an open subgroup of
  $\mathbb{G}$, and that $\mathrm{FT}_p^{\mathbb{E}} =
 \mathrm{FT}_p^{\mathbb{Z}}$\ if\ $\mathbb{E}$ is a nontrivial torsion free abelian group with the discrete topology.

 In this paper we characterize $\mathcal{FT}_p^{\mathbb{G}}$ partly as follows.
 In Section 2, we show for every locally compact abelian group
 $\mathbb{G}$ that
 $\|\,T\,|\mathcal{FT}_p^{\mathbb{R}\times\mathbb{G}}\|$ is equivalent to
 $\|\,T\,|\mathcal{FT}_p^{\mathbb{Z}\times\mathbb{G}}\|$,
  and
 $\|\,T\,|\mathcal{FT}_p^{\mathbb{Z}^n\times\mathbb{G}}\| =
 \|\,T\,|\mathcal{FT}_p^{\mathbb{Z}\times\mathbb{G}}\|$.
 The Cartesian product means the direct sum. By applying these, we easily obtain some results
 in the above paragraph and the relation
 $\mathcal{FT}_p^{\mathbb{R}^k\times\mathbb{Z}^l\times\mathbb{T}^m\times\mathbb{G}}
 =\mathcal{FT}_p^{\mathbb{Z}\times\mathbb{G}}$ for any nonnegative
 integers $k$, $l$, $m$ with $k+l+m \geq 1$.  In Section 3, we combine the results of Section 2 with
 the properties of a locally compact abelian group to
 show that
 $\mathcal{FT}_p^{\mathbb{R}^k\times\mathbb{Z}^l\times\mathbb{F}}
 = \mathcal{FT}_p^{\mathbb{Z}}$ for any compact abelian
 group $\mathbb{F}$ with finitely many components. And we show that if $\mathbb{G} \cong
 \mathbb{R}^k
\times$[torsion-free group with discrete topology]$\times$[compact
group with finite components] then $\mathcal{FT}_p^{\mathbb{G}}
=\mathcal{FT}_p^{\mathbb{Z}}$. If the factor group of a locally
compact abelian group $\mathbb{G}$ with respect to the component
of the identity element is in the form of [torsion-free
group]$\times$[finite group] with discrete topology, then we have
$\mathcal{FT}_p^{\mathbb{G}} =\mathcal{FT}_p^{\mathbb{Z}}$.

From now on, $ X $ and $ Y $ are Banach spaces and $ T : X \to Y $
is a bounded linear operator. We denote the dual group of
$\mathbb{G}$ by $\mathbb{G}^{'}$. But we use the fact that
$\mathbb{R}^{'}=\mathbb{R}$ and $\mathbb{Z}^{'}=\mathbb{T}$. We
use the abbreviation LCA for `locally compact abelian'. The term
`isomorphic' means `topologically and algebraically isomorphic'.
The integral for a vector valued function is the Bochner integral.

\section{Classifying $\mathcal{FT}_p^{\mathbb{G}}$ via direct sum.}
First we observe that the proof in \cite{Gar} of the fact
$\mathrm{FT}_p^{\mathbb{R}}=\mathrm{FT}_p^{\mathbb{Z}}$ can be
modified to yield the following:
\begin{prop} \label{EQRZ} For any LCA group $\mathbb{G}$, we have
the inequalities
$$\|\,T\,|\mathcal{FT}_p^{\mathbb{R}\times\mathbb{G}}\,\|\, \leq
\|\,T\,|\mathcal{FT}_p^{\mathbb{Z}\times\mathbb{G}}\,\|\, \leq \,
\frac\pi2\,
\|\,T\,|\mathcal{FT}_p^{\mathbb{R}\times\mathbb{G}}\,\|,$$ and
hence
$$\mathcal{FT}_{p}^{\mathbb{R} \times \mathbb{G}} =
\mathcal{FT}_{p}^{\mathbb{Z} \times \mathbb{G }}.$$
\end{prop}

\begin{proof}
For arbitrary $\delta > 0 $, let $$ {\mathbf{f}}(s,t) = \sum_{m}
\chi_{[\delta (m-\frac{1}{2})\,,\, \delta (m+\frac{1}{2}) ]}(s)\,
{\mathbf{g}}_{m} (t),$$ where the summation is over $\mathbb{Z}$,
${\mathbf{g}}_{m}$ is an $X$-valued simple function on
$\mathbb{G}$ and ${\mathbf{g}}_{m}=0$ except for finitely many
$m$. Note that the set of all $X$-valued functions of the form of
$\mathbf{f}$ is dense in $[L_p(\mathbb{R}\times\mathbb{G}),X]$. We
compute
\begin{equation}
\label{lpnorm}
\begin{split}
\|\mathbf{f}| L_{p} (\mathbb{R} \times \mathbb{G})\|
&=\Big{(}\int_{\mathbb{R} \times \mathbb{G}}
\|\sum_{m}\chi_{[\delta (m-\frac12)\,,\, \delta (m+\frac{1}{2})
]}(s)\, \mathbf{g}_m
(t)\|^{p} \,ds\,dt\Big{)}^{\frac{1}{p}} \\
&=\delta^{\frac{1}{p}}\Big{(}\int_{\mathbb{G}}\,\sum_{m}
\|{\mathbf{g}}_{m} (t)\|^{p} \,dt\Big{)}^{\frac{1}{p}}.\\
\end{split}
\end{equation}
Since $\int_{\delta(m-1/2)}^{\delta(m+1/2)} \exp(i\tilde{s}s) \,ds
\,=\,\exp(im\delta\tilde{s}) \,\frac{\sin(\frac{\delta}{2}
\tilde{s})}{\tilde{s}/2}$, we have
\begin{equation*}
\begin{split}
\widehat{T\mathbf{f}}(\tilde{s}\,,\,\tilde{t}) &=
\int_{\mathbb{G}}\,\int_{\mathbb{R}}\,\sum_{m} T\mathbf{g}_{m} (t)
\, \chi_{[\delta (m-\frac{1}{2})\,,\, \delta (m +\frac{1}{2})
]}(s)
\exp(i\tilde{s}s)\,(t,\tilde{t})\,ds\,dt \\
&=
\int_{\mathbb{G}}\,\sum_{m}\,T\mathbf{g}_{m}(t)\,(t,\tilde{t})\,\exp(i\delta
m\tilde{s})\,\frac{\sin(\frac{\delta}{2}
\tilde{s})}{\frac12\tilde{s}} \,dt. \\
\end{split}
\end{equation*}
Hence
\begin{equation*}
\begin{split}
\|[F^{\mathbb{R} \times \mathbb{G}} \,&,T]\mathbf{f}\|^{p^{'}} =
\int_{\mathbb{G}^{'}}\,\int_{\mathbb{R}}\,\Big{\|}\int_{\mathbb{G}}\,\sum_{m}\,T{\mathbf{g}}_{m}(t)
\,(t,\tilde{t})\,\exp(i\delta m \tilde{s})
\,\frac{\sin(\frac{\delta}{2} \tilde{s})}{\frac12\tilde{s}}
\,dt\,\Big{\|}^{p^{'}}\,\frac1{2\pi}d\tilde{s}\,d\tilde{t} \\ &=
\delta^{p^{'}}\,\int_{\mathbb{G}^{'}}\,\int_{\mathbb{R}}\,\Big{\|}\int_{\mathbb{G}}\,\sum_{m}\,T\mathbf{g}_{m}(t)
\,(t,\tilde{t})\,\exp(i\delta
m\tilde{s})\,\frac{\sin(\frac{\delta}{2}\tilde{s})}{\frac{\delta}{2}
\tilde{s}}\,dt\,\Big{\|}^{p^{'}}\,\frac1{2\pi}d\tilde{s}\,d\tilde{t}
\\ &=
\delta^{p^{'}-1}\,\int_{\mathbb{G}^{'}}\,\int_{\mathbb{R}}\,\Big{|}\frac{\sin(\frac12\tilde{s})}{\frac12\tilde{s}}\Big{|}^{p^{'}}
\Big{\|}\int_{\mathbb{G}}\,\sum_{k}\,(\sum_{m}\,T\mathbf{g}_{m}(t)\,\chi_{\{m\}}(k))\,(t,\tilde{t})\,\,exp(ik\tilde{s})
\,dt\Big{\|}\,\frac1{2\pi}d\tilde{s}\,d\tilde{t} \\ &=
\,\delta^{p^{'}-1}\,\int_{\mathbb{G}^{'}}\,\int_{-\pi}^{\pi}
\,\sum_{n}
\,\Big{|}\frac{\sin(\frac12\tilde{s})}{\frac12\tilde{s}-n\pi}\Big{|}^{p^{'}}
\big{\|}[F^{\mathbb{Z}\times\mathbb{G}}\,,T](\sum_{m}\,{\mathbf{g}}_{m}\,\chi_{\{m\}})\big{\|}^{p^{'}}\,\frac1{2\pi}d\tilde{s}\,d\tilde{t}
\\ &\leq \delta^{p^{'}-1}
\|\,T\,|\mathcal{FT}_p^{\mathbb{Z} \times \mathbb{G}}\|^{p^{'}}
\Big{(}\int_{\mathbb{G}} \sum_{m} \|{\mathbf{g}}_{m} (t) \|^{p}
dt\Big{)}^{\frac{p^{'}}{p}} \\
&= \|\,T\,|\mathcal{FT}_p^{\mathbb{Z} \times \mathbb{G}}\|^{p^{'}}
\,\|\, \mathbf{f}\,| L_{p}(\mathbb{R} \times \mathbb{G})
\|^{p^{'}} \qquad \mbox{by\ } (\ref{lpnorm}), \\
\end{split}
\end{equation*}
where we used the inequality  $\sum_{n}
|\frac{\sin\tilde{s}}{\tilde{s}-n\pi}|^{p^{'}} \,\leq \,1 $ for
any real $\tilde{s} \neq n\pi $ and $2\leq p^{'}< \infty$, see
\cite{Jo}.

Therefore we have $\|\,T\,|\mathcal{FT}_p^{\mathbb{R} \times
\mathbb{G}}\| \leq \|\,T\,|\mathcal{FT}_p^{\mathbb{Z} \times
\mathbb{G}}\|$ and $\mathcal{FT}_p^{\mathbb{Z}\times\mathbb{G}}
\,\subseteq\,\mathcal{FT}_p^{\mathbb{R}\times\mathbb{G}}$. \\

For the right inequality let $\mathbf{f} (k,t) = \sum_{m}
\mathbf{g}_{m}(t)\,\chi_{\{m\}}(k)$ where $\mathbf{g}_{m}$ is an
$X$-valued simple function and $\mathbf{g}_{m}=0$ except for
finitely many $m$. By density argument \ it's enough to consider
$\mathbf{f}$ of the above form. Now we have
\begin{equation}
\label{lpnorm2}
\begin{split}
\|\mathbf{f}| L_{p} (\mathbb{Z} \times
\mathbb{G})\|\,&=\,\Big{(}\int_{\mathbb{G}}\,\sum_{k}\,\|\sum_{m}\mathbf{g}_{m}(t)
\,\chi_{\{m\}} (k)\,\|^{p}\,dt\Big{)}^{\frac1p} \\ &=\,\Big{(}
\int_{\mathbb{G}}\,\sum_{m}
\|\,\mathbf{g}_{m}(t)\,\|^{p}\,dt\Big{)}^{\frac1p}; \ \mbox{\ and} \\
\end{split}
\end{equation}
\begin{equation*}
\begin{split}
\widehat{T\mathbf{f}}(\tilde{s}\,,\tilde{t})\, &=\,
\int_{\mathbb{G}}\sum_{k}\,\big{(}\sum_{m}T\mathbf{g}_{m}(t)\,\chi_{\{m\}}(k)\big{)}\exp(i\tilde{s}k)\,(t,\tilde{t})\,dt
\\ &= \,\int_{\mathbb{G}} \sum_{m}
T\mathbf{g}_{m}(t)\,\exp(i\tilde{s}m)\,(t,\tilde{t})\,dt.\\
\end{split}
\end{equation*}
We use the identity, $ \exp(im\tilde{s}) \, = \,\int_{\mathbb{R}}
\frac{\tilde{s}/2}{\sin(\tilde{s}/2)}
\,\chi_{[m-1/2\,,\,m+1/2]}(s)\,\exp(is\tilde{s})\,ds $, to obtain
the following inequality:
\begin{equation*}
\begin{split}
\|[&F^{\mathbb{Z}\times\mathbb{G}}\,,T]\mathbf{f}\|^{p^{'}}\\ &=
\, \frac1{2\pi} \!\int_{\mathbb{G}^{'}}\!\int_{-\pi}^{\pi}
\Big{\|}\!\int_{\mathbb{G}}\sum_{m}
T\mathbf{g}_{m}(t)\,(t,\tilde{t})\!\int_{\mathbb{R}}\frac{\tilde{s}/2}{\sin(\tilde{s}/2)}
\,\chi_{[m-1/2\,,\,m+1/2]}(s)\,\exp(i\tilde{s}s)\,ds\,dt\,\Big{\|}^{p^{'}}\,d\tilde{s}\,d\tilde{t}
\\
&\leq \, \frac1{2\pi}\int_{\mathbb{G}^{'}}\int_{-\pi}^{\pi}
\Big{|}\frac{\tilde{s}/2}{\sin(\tilde{s}/2)}\Big{|}^{p^{'}}
\Big{\|}\!\int_{\mathbb{G}}\!\int_{\mathbb{R}} \sum_{m}
T\mathbf{g}_{m}(t)
\,\chi_{[m-\frac12\,,\,m+\frac12]}(s)\,(t,\tilde{t})\,\exp(i\tilde{s}s)\,ds\,dt
\Big{\|}^{p^{'}}d\tilde{s}\,d\tilde{t} \\
&\ \ \ \ \ \ \ \ \hspace{1in} ( \ \mbox{because }
\Big{|}\frac{\tilde{s}/2}{\sin(\tilde{s}/2)}\Big{|} \,\leq
\,\frac\pi2\ \,\mbox{for} \ -\pi\,\leq\,\tilde{s}\,\leq\,\pi \ ) \\
&\leq
\,\frac1{2\pi}\,\Big{(}\frac\pi2\Big{)}^{p^{'}}\int_{\mathbb{G}^{'}}\,\int_{\mathbb{R}}
\Big{\|}\int_{\mathbb{G}}\int_{\mathbb{R}} \sum_{m}
T\mathbf{g}_{m}(t)
\,\chi_{[m-\frac12\,,\,m+\frac12]}(s)\,(t,\tilde{t})\,\exp(i\tilde{s}s)\,ds\,dt
\Big{\|}^{p^{'}}d\tilde{s}\,d\tilde{t} \\
&=\,\Big{(}\frac\pi2\Big{)}^{p^{'}}\,\big{\|}[F^{\mathbb{R}\times\mathbb{G}}\,,T]\big{(}\sum_{m}
\mathbf{g}_{m}\,\chi_{[m-\frac12\,,\,m+\frac12]}\big{)}\,\big{\|}^{p^{'}} \\
&\leq
\,\Big{(}\frac\pi2\Big{)}^{p^{'}}\,\|\,T\,|\mathcal{FT}_p^{\mathbb{R}\times\mathbb{G}}\|^{p^{'}}
\,\big{\|}\sum_{m}
\mathbf{g}_{m}\,\chi_{[m-\frac12\,,\,m+\frac12]}\,|L_{p}
(\mathbb{R} \times \mathbb{G})\big{\|}^{p^{'}} \\
&\ \ (\mbox{\,since } \big{\|}\sum_{m}
\mathbf{g}_{m}\,\chi_{[m-\frac12\,,\,m+\frac12]}\,|L_p(\mathbb{R}\times\mathbb{G})\big{\|}\,=\,
\Big{(}\int_{\mathbb{G}}\sum_{m}\|\mathbf{g}_{m}(t)\,\|^p\,dt\Big{)}^{\frac1p})\\
&=\,\Big{(}\frac\pi2\Big{)}^{p^{'}}\,\|\,T\,|\mathcal{FT}_p^{\mathbb{R}\times\mathbb{G}}\|^{p^{'}}
\,\|\mathbf{f}\,|L_{p} (\mathbb{Z} \times \mathbb{G})\|^{p^{'}}
\qquad \ \mbox{by}\ (\ref{lpnorm2}). \\
\end{split}
\end{equation*}
Therefore we have
$\|\,T\,|\mathcal{FT}_p^{\mathbb{Z}\times\mathbb{G}}\| \leq \,
\frac\pi2\, \|\,T\,|\mathcal{FT}_p^{\mathbb{R}\times\mathbb{G}}\|$
 and $\mathcal{FT}_p^{\mathbb{R}\times\mathbb{G}}
\,\subseteq\,\mathcal{FT}_p^{\mathbb{Z}\times\mathbb{G}}$. \\
This completes the proof of Proposition \ref{EQRZ}.
\end{proof}
According to Theorem 6.3 of \cite{Gar}, for any LCA group
$\mathbb{G}$ and operator $T$,
$\|\,T\,|\mathcal{FT}_p^{\mathbb{G}}\|=\|\,T^{'}\,|\mathcal{FT}_p^{\mathbb{G}^{'}}\|$.
By applying this property, we have the following:
\begin{prop}
\label{EQRT} For every LCA group $\mathbb{G}$,
$\mathcal{FT}_{p}^{\mathbb{R} \times \mathbb{G}} =
\mathcal{FT}_{p}^{\mathbb{T} \times \mathbb{G }}$.
\end{prop}

\begin{proof}

By Proposition \ref{EQRZ} we have
$$
\|\,T\,|\mathcal{FT}_p^{\mathbb{R}\times\mathbb{G}}\|\,=\,\|\,T^{'}\,|\mathcal{FT}_p^{\mathbb{R}\times\mathbb{G}^{'}}\|\,\leq
\,\|\,T^{'}\,|\mathcal{FT}_p^{\mathbb{Z}\times\mathbb{G}^{'}}\|\,
= \,\|\,T\,|\mathcal{FT}_p^{\mathbb{T}\times\mathbb{G}}\|.$$
Similarly
$$\|\,T\,|\mathcal{FT}_p^{\mathbb{T}\times\mathbb{G}}\| \,\leq
\,\frac\pi2\,\|\,T\,|\mathcal{FT}_p^{\mathbb{R}\times\mathbb{G}}\|.$$
\end{proof}

From Proposition \ref{EQRZ} and \ref{EQRT}, we conclude that
$\mathcal{FT}_{p}^{\mathbb{R} \times \mathbb{G}} =
\mathcal{FT}_{p}^{\mathbb{Z} \times \mathbb{G}} =
\mathcal{FT}_{p}^{\mathbb{T} \times \mathbb{G}} $ for every LCA
group $\mathbb{G}$. And we have the following two corollaries.

\begin{cor}
\label{NRTZ} $\mathcal{FT}_{p}^{{\mathbb{R}}^{n}}
\,=\,\mathcal{FT}_{p}^{{\mathbb{Z}}^{n}}\,=\,\mathcal{FT}_{p}^{{\mathbb{T}}^{n}}$\
for every integer $n\geq1$.
\end{cor}
\begin{proof}
By Proposition \ref{EQRZ} we have
$$\mathcal{FT}_{p}^{{\mathbb{R}}^{n}}\,=\,\mathcal{FT}_{p}^{\mathbb{R}\times{\mathbb{R}}^{n-1}}\,=\,
\mathcal{FT}_{p}^{\mathbb{Z}\times{\mathbb{R}}^{n-1}}.  $$ We continue this to obtain \\
$$\mathcal{FT}_{p}^{\mathbb{R}\times\mathbb{Z}\times{\mathbb{R}}^{n-2}}\,=\,
\mathcal{FT}_{p}^{\mathbb{Z}\times\mathbb{Z}\times{\mathbb{R}}^{n-2}}\,=\,...\,=\,
\mathcal{FT}_{p}^{{\mathbb{Z}}^{n}}.$$ Similarly we apply
Proposition \ref{EQRT} \ $n$ \ times to obtain
$$\mathcal{FT}_{p}^{{\mathbb{R}}^{n}}
\,=\,\mathcal{FT}_{p}^{{\mathbb{T}}^{n}}.
$$
\end{proof}
\begin{cor}
\label{NRTZG} $\mathcal{FT}_{p}^{{\mathbb{R}}^{n}\times
\mathbb{G}} \,=\,\mathcal{FT}_{p}^{{\mathbb{Z}}^{n}\times
\mathbb{G}}\,=\,\mathcal{FT}_{p}^{{\mathbb{T}}^{n}\times
\mathbb{G}}$ \ for every LCA group $\mathbb{G}$.
\end{cor}
\begin{proof}
The proof is similar to that of Corollary \ref{NRTZ}.
\end{proof}

\begin{lem}
\label{ZNG} For every LCA group $\mathbb{G}$ and every positive
integer $n$,
$$\|\,T\,|\mathcal{FT}_p^{\mathbb{Z}^n\times\mathbb{G}}\|
\,=\,\|\,T\,|\mathcal{FT}_p^{\mathbb{Z}\times\mathbb{G}}\|$$ and
therefore $\mathcal{FT}_{p}^{{\mathbb{Z}}^{n}\times
\mathbb{G}}\,=\,\mathcal{FT}_{p}^{\mathbb{Z}\times\mathbb{G}}$.
\end{lem}
\begin{proof}
Let
$\mathbf{f}(k,t)\,=\,\sum_{m}\mathbf{g}_{m}(t)\,\chi_{\{m\}}(k)$
\,,where $\mathbf{g}_{m}$  is an $X$-valued simple function on
$\mathbb{G}$ and the summation is on a finite subset of
$\mathbb{Z}$, then $$\|\mathbf{f}| L_{p} (\mathbb{Z} \times
\mathbb{G})\|\,=\,\Big{(} \int_{\mathbb{G}}\,\sum_{m}
\|\,\mathbf{g}_{m}(t)\,\|^{p}\,dt\Big{)}^{\frac1p}.$$ We have
$$[F^{\mathbb{Z}\times\mathbb{G}}\,,\,T]\mathbf{f}(\tilde{s},\tilde{t})
\,=\,\int_{\mathbb{G}}\sum_{k}
(\sum_{m}T\mathbf{g}_{m}(t)\,\chi_{\{m\}}(k))\,\exp(i\tilde{s}k)\,(t,\tilde{t})\,dt$$
 and
\begin{equation*}
\begin{split}
\|[&F^{\mathbb{Z}\times\mathbb{G}}\,,\,T]\mathbf{f}\|^{p^{'}}
\,=\int_{\mathbb{G}^{'}}\int_{\mathbb{T}}
\Big{\|}\int_{\mathbb{G}}\sum_{k}
(\sum_{m}T\mathbf{g}_{m}(t)\,\chi_{\{m\}}(k))\,\exp(i\tilde{s}k)\,(t,\tilde{t})
\,dt\Big{\|}^{p{'}}\,d\tilde{s}\,d\tilde{t} \\
&=\int_{\mathbb{G}^{'}}\!\int_{\mathbb{T}}\!\int_{\mathbb{T}}
\Big{\|}\!\int_{\mathbb{G}}\sum_{k,l}\big{(}\sum_{m}T\mathbf{g}_m(t)\,\chi_{\{m\}}(k)\chi_{\{0\}}(l)\big{)}
\,\exp(i\tilde{s}k)\,\exp(i\xi
l)\,(t,\tilde{t})\,dt\Big{\|}^{p^{'}}\,d\tilde{s}\,d\xi\,d\tilde{t} \\
&\leq\,\|\,T\,|\mathcal{FT}_p^{\mathbb{Z}^2\times\mathbb{G}}\|^{p^{'}}\Big{(}\int_{\mathbb{G}}\sum_{k,l}
\big{\|}\sum_m\mathbf{g}_m(t)\,\chi_{\{m\}}(k)\chi_{\{0\}}(l)\big{\|}^pdt\Big{)}^{\frac{p^{'}}p}
\\ &=\,\|\,T\,|\mathcal{FT}_p^{\mathbb{Z}^2\times\mathbb{G}}\|^{p^{'}}\Big{(} \int_{\mathbb{G}}\,\sum_{m}
\|\,\mathbf{g}_{m}(t)\,\|^{p}\,dt\Big{)}^{\frac{p^{'}}p},\\
\end{split}
\end{equation*}
so
$\|\,T\,|\mathcal{FT}_p^{\mathbb{Z}\times\mathbb{G}}\|\,\leq\,\|\,T\,|\mathcal{FT}_p^{\mathbb{Z}^2\times\mathbb{G}}\|$.

Conversely let $$\mathbf{f}(k_1,k_2,t)\,=
\,\sum_{l_1,l_2}\mathbf{g}_{l_1,l_2}(t)\,\chi_{\{l_1\}}(k_1)\,\chi_{\{l_2\}}(k_2),$$
where the summation is on a finite subset of
$\mathbb{Z}\times\mathbb{Z}$, then we have
\begin{equation}
\label{lpz2norm}
\|\mathbf{f}|L_p(\mathbb{Z}^2\times\mathbb{G})\|\,=\,\Big{(}\int_{\mathbb{G}}\sum_{l_1,l_2}
\|\mathbf{g}_{l_1,l_2}(t)\|^pdt\Big{)}^{\frac1p}.\end{equation}
Since
$[F^{\mathbb{Z}^2\times\mathbb{G}}\,,\,T]\mathbf{f}(\tilde{s}_1,\tilde{s}_2,\tilde{t})\,=\,
\int_{\mathbb{G}}\sum_{l_1,l_2}T\mathbf{g}_{l_1,l_2}(t)\,\exp(i\tilde{s}_1l_1)\,\exp(i\tilde{s}_2l_2)
\,(t,\tilde{t})\,dt$ we have
\begin{equation}
\label{LEQNZG}
\begin{split}
\|[F^{\mathbb{Z}^2\times\mathbb{G}}\,&,\,T]\mathbf{f}\|^{p^{'}}=
\\
&\ \ \ \ \
\int_{\mathbb{G}^{'}}\int_{\mathbb{T}}\int_{\mathbb{T}}\Big{\|}\!\int_{\mathbb{G}}\sum_{l_1,l_2}
T\mathbf{g}_{l_1,l_2}(t)\,\exp(i\tilde{s}_1l_1)\,\exp(i\tilde{s}_2l_2)
\,(t,\tilde{t})\,dt\Big{\|}^{p^{'}}d\tilde{s}_1\,d\tilde{s}_2\,d\tilde{t}
\\
&\ \ \ \  \ \ \ \ \ \ \ \ (\ \mbox{we let } A \,=\,2\:max\{l_1\}\,+1 \ )\\
&=\!\int_{\mathbb{G}^{'}}\!\int_{\mathbb{T}}\!\int_{\mathbb{T}}\Big{\|}\!\int_{\mathbb{G}}\sum_{l_1,l_2}
T\mathbf{g}_{l_1,l_2}(t)\,\exp(i\tilde{s}_1l_1)\,\exp(i\tilde{s}_2Al_2)
\,(t,\tilde{t})\,dt\Big{\|}^{p^{'}}d\tilde{s}_1\,d\tilde{s}_2\,d\tilde{t}
\\
&=\!\int_{\mathbb{G}^{'}}\!\int_{\mathbb{T}}\!\int_{\mathbb{T}}\Big{\|}\!\int_{\mathbb{G}}\sum_{l_1,l_2}
T\mathbf{g}_{l_1,l_2}(t)\,\exp(i\tilde{s}_1l_1)\,\exp(i(\tilde{s}_1+\tilde{s}_2)Al_2)
\,(t,\tilde{t})\,dt\Big{\|}^{p^{'}}d\tilde{s}_1\,d\tilde{s}_2\,d\tilde{t}
\\
&=\!\int_{\mathbb{G}^{'}}\!\int_{\mathbb{T}}\!\int_{\mathbb{T}}\Big{\|}\!\int_{\mathbb{G}}\sum_{l_1,l_2}
T\mathbf{g}_{l_1,l_2}(t)\,\exp(i\tilde{s}_1(l_1+Al_2))\,\exp(i\tilde{s}_2Al_2)
\,(t,\tilde{t})\,dt\Big{\|}^{p^{'}}d\tilde{s}_1\,d\tilde{s}_2\,d\tilde{t}
\\
&=\!\int_{\mathbb{T}}\!\int_{\mathbb{G}^{'}}\!\int_{\mathbb{T}}\Big{\|}\!\int_{\mathbb{G}}\sum_{l_1,l_2}
T(\mathbf{g}_{l_1,l_2}(t)\,\exp(i\tilde{s}_2Al_2))\,\exp(i\tilde{s}_1(l_1+Al_2))
\,(t,\tilde{t})\,dt\Big{\|}^{p^{'}}d\tilde{s}_1\,d\tilde{t}\,d\tilde{s}_2
\\
&=\int_{\mathbb{T}}\big{\|}[F^{\mathbb{Z}\times\mathbb{G}}\,,\,T]\big{(}\sum_{l_1,l_2}
(\mathbf{g}_{l_1,l_2}(t)\,\exp(i\tilde{s}_2Al_2)\big{)}\,\chi_{\{l_1+Al_2\}}\big{\|}^{p^{'}}d\tilde{s}_2\\
&\leq
\int_{\mathbb{T}}\Big{[}\|\,T\,|\mathcal{FT}_p^{\mathbb{Z}\times\mathbb{G}}\|^{p^{'}}\Big{(}\int_{\mathbb{G}}\sum_{l_1+Al_2}
\|\mathbf{g}_{l_1,l_2}(t)\,\exp(i\tilde{s}_2Al_2)\|^pdt\Big{)}^{p^{'}/p}\Big{]}\,d\tilde{s}_2\\
&=\|\,T\,|\mathcal{FT}_p^{\mathbb{Z}\times\mathbb{G}}\|^{p^{'}}\big{(}\int_{\mathbb{G}}\sum_{l_1+Al_2}
\|\mathbf{g}_{l_1,l_2}(t)\|^pdt\big{)}^{p^{'}/p} \\
&=\|\,T\,|\mathcal{FT}_p^{\mathbb{Z}\times\mathbb{G}}\|^{p^{'}}\big{(}\int_{\mathbb{G}}\sum_{l_1,l_2}
\|\mathbf{g}_{l_1,l_2}(t)\|^pdt\big{)}^{p^{'}/p} \\
&=\|\,T\,|\mathcal{FT}_p^{\mathbb{Z}\times\mathbb{G}}\|^{p^{'}}\|\mathbf{f}|L_p(\mathbb{Z}^2\times\mathbb{G})\|^{p^{'}}
\qquad \mbox{by (\ref{lpz2norm})}, \\
\end{split}
\end{equation}
so it follows that
$\|\,T\,|\mathcal{FT}_p^{\mathbb{Z}^2\times\mathbb{G}}\|\,\leq\,\|\,T\,|\mathcal{FT}_p^{\mathbb{Z}\times\mathbb{G}}\|$.

Therefore we have
$\|\,T\,|\mathcal{FT}_p^{\mathbb{Z}^2\times\mathbb{G}}\|\,=\,\|\,T\,|\mathcal{FT}_p^{\mathbb{Z}\times\mathbb{G}}\|
\,.$ A simple generalization yields, for all positive integers $n$
that
\begin{equation}
\label{NZEQ}
\|\,T\,|\mathcal{FT}_p^{\mathbb{Z}^n\times\mathbb{G}}\|\,=\,\|\,T\,|\mathcal{FT}_p^{\mathbb{Z}\times\mathbb{G}}\|
\end{equation}
hence that $
\mathcal{FT}_p^{\mathbb{Z}^n\times\mathbb{G}}\,=\,\mathcal{FT}_p^{\mathbb{Z}\times\mathbb{G}}
$.
\end{proof}
\begin{rem}
In fact, a dissipative group $\mathbb{A}$ satisfies the property
that $\|\,T\,|\mathcal{FT}_p^{{\mathbb{A}}^n\times\mathbb{G}}\|\,=
\,\|\,T\,|\mathcal{FT}_p^{\mathbb{A}\times\mathbb{G}}\|  $ \ for
any positive integer $n$ and for any LCA group $\mathbb{G}$. The
definition and properties of a dissipative group are introduced in
\cite{Gar}. The idea of the proof of Lemma \ref{ZNG} goes back to
\cite{Gar}.
\end{rem}

\begin{cor}
$\|\,T\,|\mathcal{FT}_p^{\mathbb{T}^n\times\mathbb{G}}\|
\,=\,\|\,T\,|\mathcal{FT}_p^{\mathbb{T}\times\mathbb{G}}\|$ and so
$\mathcal{FT}_{p}^{{\mathbb{T}}^{n}\times
\mathbb{G}}\,=\,\mathcal{FT}_{p}^{\mathbb{T}\times\mathbb{G}}$.
\end{cor}
\begin{proof}
This follows from Lemma \ref{ZNG} and a duality argument.
\end{proof}
\begin{prop}
\label{HNZ}
 Let $\mathbb{H}$ be $\mathbb{R}$, $\mathbb{Z}$ or
$\mathbb{T}$, and $\mathbb{G}$ be an LCA group. Then
$\mathcal{FT}_p^{\mathbb{H}^n\times \mathbb{G}
}\,=\,\mathcal{FT}_p^{\mathbb{Z}\times\mathbb{G}}$ for any integer
$n\, \geq\,1$. In particular,
$\mathcal{FT}_p^{\mathbb{H}^n}\,=\,\mathcal{FT}_p^{\mathbb{Z}}$.
\end{prop}
\begin{proof}
It follows from Corollary $\ref{NRTZG}$ and Lemma $\ref{ZNG}$ that
$\mathcal{FT}_p^{\mathbb{H}^n\times \mathbb{G}
}\,=\,\mathcal{FT}_p^{\mathbb{Z}\times\mathbb{G}}$. And if
$\mathbb{G}$ is the trivial group then we have
$\mathcal{FT}_p^{\mathbb{H}^n}\,=\,\mathcal{FT}_p^{\mathbb{Z}}$.
\end{proof}

Notation: From now on, when $\mathbb{H}$ is $\mathbb{R}$,
$\mathbb{Z}$ or $\mathbb{T}$ we denote
$\mathcal{FT}_p^{\mathbb{H}}$ by $\mathcal{FT}_p$.

\begin{cor}
$\mathcal{FT}_p^{\mathbb{R}^a\times\mathbb{Z}^b\times\mathbb{T}^c}=\,\mathcal{FT}_p
$ for any nonnegative integers $a\,,\,b\,,\,c$ with
 $a+b+c\,\geq\,1$.
\end{cor}
\begin{proof}
Without loss of generality we assume $a, b, c > 0$. By applying
Proposition $\ref{HNZ}$, we have
$$\mathcal{FT}_p^{\mathbb{R}^a\times\mathbb{Z}^b\times\mathbb{T}^c}
=\,\mathcal{FT}_p^{\mathbb{Z}^{b+1}\times\mathbb{T}^c}=\,\mathcal{FT}_p^{\mathbb{Z}\times\mathbb{T}^c}
=\,\mathcal{FT}_p^{\mathbb{Z}^2}=\,\mathcal{FT}_p^{\mathbb{Z}}.$$
\end{proof}
\begin{thm}
\label{HNG} Let $\mathbb{G}$ be an LCA group then
$\mathcal{FT}_p^{\mathbb{R}^a\times\mathbb{Z}^b\times\mathbb{T}^c\times\mathbb{G}}=\,
\mathcal{FT}_p^{\mathbb{Z}\times\mathbb{G}}$ for any nonnegative
integers $a\,,\,b\,,\,c$ with $a+b+c\,\geq\,1$.
\end{thm}
\begin{proof}
By Proposition $\ref{HNZ}$ we have
$$
\mathcal{FT}_p^{\mathbb{R}^a\times\mathbb{Z}^b\times\mathbb{T}^c\times\mathbb{G}}=\,
\mathcal{FT}_p^{\mathbb{Z}^{1+b}\times\mathbb{T}^c\times\mathbb{G}}=
\,\mathcal{FT}_p^{\mathbb{Z}^{2+b}\times\mathbb{G}}=\,\mathcal{FT}_p^{\mathbb{Z}\times\mathbb{G}}.$$
\end{proof}
\section{Search for $\mathbb{G}$ satisfying $\mathcal{FT}_p^{\mathbb{G}} = \mathcal{FT}_p^{\mathbb{Z}}$
and Further classification.}
We notice that to solve the problem of deciding whether
$\mathcal{FT}_p^{\mathbb{G}}$ depends on $\mathbb{G}$, we should
solve that problem on the compact abelian groups since any
 compactly generated LCA group is a product of finite number of
$\mathbb{R}$'s and $\mathbb{Z}$'s and a compact group \cite{Hw}.
The dual group of a compact group has the discrete topology
\cite{Rud}. More over if $\mathbb{G}$ is a compact abelian group
then $\mathbb{G}$ is connected iff $\mathbb{G}^{'}$ is
tortion-free \cite{Hw}. This fact gives a clue to the results in
this section.

First, we introduce Weil's formula which factorizes an integration
on $\mathbb{G}$ into double ones on a closed subgroup and its
factor group.
\begin{thm}\cite{Re}
\label{Weil}
Let $\mathbb{G}$ be an LCA group and $\mathbb{H}$ a
closed subgroup. Then there are Haar measures $\mu_{\mathbb{G}}$,
$\mu_{\mathbb{H}}$ and $\mu_{\mathbb{G}/\mathbb{H}}$ such that
\begin{equation}
\label{Wl} \int_{\mathbb{G}}\mathbf{f}(s)\,d_{\mu_{\mathbb{G}}}(s)
\,=\int_{\mathbb{G}/\mathbb{H}}\big{(}\int_{\mathbb{H}}\mathbf{f}(s+h)\,
d_{\mu_{\mathbb{H}}}(h)\big{)}d_{\mu_{\mathbb{G}/\mathbb{H}}}(s+\mathbb{H})
\end{equation}
if $\mathbf{f}$ is a compactly supported continuous Banach
space-valued function or a nonnegative lower semi continuous
function on $\mathbb{G}$.
\end{thm}
In Theorem \ref{Weil} if any two of $\mu_{\mathbb{G}}$,
$\mu_{\mathbb{H}}$ and $\mu_{\mathbb{G}/\mathbb{H}}$ are given
then the third is determined so that the statement holds.

Andersson \cite{An} obtained the inequality,
$\|\,I_{X}\,|\mathcal{FT}_p^{\mathbb{H}}\|\,\leq\,\|\,I_{X}\,|\mathcal{FT}_p^{\mathbb{G}}\|$,
where $\mathbb{H}$ is an open subgroup of an LCA group
$\mathbb{G}$ and $I_{X}$ is the identity operator on a Banach
space $X$. By replacing $I_{X}$ with a bounded linear operator
$T:X\rightarrow Y$ in the proof of \cite{An}, we have the
following:
\begin{prop}
\label{OPP}
Let $\mathbb{H}$ be an open subgroup of an LCA group
$\mathbb{G}$ then we have
\begin{equation}
\label{OS}
\|\,T\,|\mathcal{FT}_p^{\mathbb{H}}\|\,\leq\,\|\,T\,|\mathcal{FT}_p^{\mathbb{G}}\|.
\end{equation}
\end{prop}
\vspace{0.85cm}
Now we consider torsion-free LCA groups with the
discrete topology.
\begin{cor}
For any non-trivial torsion-free LCA group $\mathbb{E}$ with the
discrete topology,
\begin{equation} \label{FZ}
\|\,T\,|\mathcal{FT}_p^{\mathbb{E}}\|=
\|\,T\,|\mathcal{FT}_p^{\mathbb{Z}}\|.
\end{equation}
\end{cor}
\begin{proof}
The proof is similar to that of
$\mathrm{FT}_p^{\mathbb{E}}=\mathrm{FT}_p^{\mathbb{Z}}$ in
\cite{An}: Since $\mathbb{Z}$ is isomorphic with an open subgroup
of $\mathbb{E}$, by applying ($\ref{OS}$) we have
$\|\,T\,|\mathcal{FT}_p^{\mathbb{Z}}\|\,\leq\,\|\,T\,|\mathcal{FT}_p^{\mathbb{E}}\|$.
On the other hand for any $X$-valued simple function $\mathbf{f}$
defined on $\mathbb{E}$ which has finite $L_p$-norm, $\mathbf{f}$
is non-zero only on a subset of an open subgroup which is
isomorphic with $\mathbb{Z}^k$ for some positive integer $k$.
Therefore we have
\begin{equation*}
\frac{\|[F^{\mathbb{E}},T]\mathbf{f}\|_{p^{'}}}{\|\mathbf{f}|L_p(\mathbb{E})\|}=
\frac{\|[F^{\mathbb{Z}^k},T]\mathbf{f}\|_{p^{'}}}{\|\mathbf{f}|L_p(\mathbb{Z}^k)\|}\leq\,
\|\,T\,|\mathcal{FT}_p^{\mathbb{Z}^{k}}\|
\end{equation*} and hence
\begin{equation*}
\|\,T\,|\mathcal{FT}_p^{\mathbb{E}}\|\,\leq
\,\sup_{k}\|\,T\,|\mathcal{FT}_p^{\mathbb{Z}^{k}}\|\,=\,\|\,T\,|\mathcal{FT}_p^{\mathbb{Z}}\|.
\end{equation*}
Thus the equality (\ref{FZ}) follows.
\end{proof}
\ \\
\begin{lem}
\label{FCC} Let $\mathbb{H}$ be a closed subgroup of an LCA group
$\mathbb{G}$ such that $\mathbb{G}/\mathbb{H}$ is finite. If $n$
is the cardinal number of $\mathbb{G}/\mathbb{H}$ then
\begin{equation}
\label{FC}
\|\,T\,|\mathcal{FT}_p^{\mathbb{H}}\|\,\leq\,\|\,T\,|\mathcal{FT}_p^{\mathbb{G}}\|\,\leq\,
\|\,T\,|\mathcal{FT}_p^{\mathbb{H}}\|\cdot n^{1/{p^{'}}},
\end{equation}
hence $\mathcal{FT}_p^{\mathbb{G}} = \mathcal{FT}_p^{\mathbb{H}}$.
\end{lem}
\begin{proof}
Since $\mathbb{G}/\mathbb{H}$ is finite, $\mathbb{H}$ is open and
by the inequlity (\ref{OS}) we have the left side. The proof of
the right side is as follows. Let us choose $s_1\,,s_2\,,...,s_n$
in $\mathbb{G}$ such that
$\mathbb{G}/\mathbb{H}=\{s_1+\mathbb{H}\,,s_2+\mathbb{H}\,,...,s_n+\mathbb{H}\}$.
The measure of $\mathbb{G}/\mathbb{H}$ is the Haar measure of unit
mass. By applying Weil's formula, for any compactly supported
continuous $X$-valued function $\mathbf{f}$ on $\mathbb{G}$ we
have
\begin{equation}
\begin{split}
\int_{\mathbb{G}}\|\mathbf{f}(s)\|^pd_{\mu_{\mathbb{G}}}(s)&=\int_{\mathbb{G}/\mathbb{H}}\Big{(}
\int_{\mathbb{H}}
\|\mathbf{f}(s+h)\|^p\,d_{\mu_{\mathbb{H}}}(h)\Big{)}\,d_{\mu_{\mathbb{G}/\mathbb{H}}}(s+\mathbb{H})
\\ &=\sum_{i=1}^n \frac1n \int_{\mathbb{H}}
\|\mathbf{f}(s_i+h)\|^p\,d_{\mu_{\mathbb{H}}}(h) \\
\end{split}
\end{equation} and
\begin{equation}
\begin{split}
[F^{\mathbb{G}},T]\mathbf{f}(\sigma)&=\sum_{i=1}^n \frac1n
\Big{(}\int_{\mathbb{H}}T\mathbf{f}(s_i+h)\,\sigma(s_i+h)\,d_{\mu_{\mathbb{H}}}(h)\Big{)} \\
&=\sum_{i=1}^n \frac{\sigma(s_i)}n
\Big{(}\int_{\mathbb{H}}T\mathbf{f}(s_i+h)\,\sigma(h)\,d_{\mu_{\mathbb{H}}}(h)\Big{)}
\end{split}
\end{equation} for $\sigma\in\mathbb{G}^{'}$. The dual group of
$\mathbb{G}/\mathbb{H}$ is isomorphic with the closed subgroup
$\mathbb{H}^{\bot}=\{\chi\in\mathbb{G}^{'}|\,\chi(h)=1 \mbox{ for
all }  h \in \mathbb{H}\}$ and $\mathbb{G}^{'}/\mathbb{H}^{\bot}$
is isomorphic with $\mathbb{H}^{'}$. Here the cardinal number of
$\mathbb{H}^{\bot}$ is $n$, the measure of $\mathbb{H}^{\bot}$ is
the counting measure which make Parseval's identity clean with
constant 1 and we write $\mathbb{H}^{\bot}$ as
$\{\tilde{\eta}_1\,,\tilde{\eta}_2\,,...,\tilde{\eta}_n\}$.
$[F^{\mathbb{G}},T](\mathbf{f})$ belongs to
$[C_0(\mathbb{G}^{'}),Y]$, so $Y$-norm of
$[F^{\mathbb{G}},T](\mathbf{f})$ is continuous and by Weil's
formula we have
\begin{equation*}
\begin{split}
\|[F^{\mathbb{G}},T]\mathbf{f}\|&=\Big{(}\int_{\mathbb{H}^{'}}\sum_{j=1}^n\big{\|}\frac1n\sum_{i=1}^n\sigma(s_i)\,\eta_j(s_i)
\int_{\mathbb{H}}T\mathbf{f}(s_i+h)\sigma(h)\,d_{\mu_{\mathbb{H}}}(h)\,\big{\|}^{p^{'}}d_{{\mu}_{\mathbb{H}^{'}}}(\sigma+\mathbb{H}^{\bot})
\Big{)}^{1/{p^{'}}}
\\&\leq\Big{[}\Big{(}\int_{\mathbb{H}^{'}}\sum_{j=1}^n \big{(}\sum_{i=1}^n \frac1n \Big{\|}\int_{\mathbb{H}}T\mathbf{f}(s_i+h)\sigma(h)\,d_{\mu_{\mathbb{H}}}(h)
\Big{\|}^{p}\big{)}^{{p^{'}}/{p}}d_{{\mu}_{\mathbb{H}^{'}}}(\sigma+\mathbb{H}^{\bot})\Big{)}^{p/{p^{'}}}\Big{]}^{1/p}
\\&\leq \Big{[}\sum_{i=1}^n \frac1n \Big{(}\int_{\mathbb{H}^{'}}\sum_{j=1}^n\Big{\|}\int_{\mathbb{H}}T\mathbf{f}(s_i+h)\sigma(h)\,
d_{\mu_{\mathbb{H}}}(h)\Big{\|}^{p^{'}}d_{{\mu}_{\mathbb{H}^{'}}}(\sigma+\mathbb{H}^{\bot})\Big{)}^{p/{p^{'}}}\Big{]}^{1/p}\\
&= \Big{[}\sum_{i=1}^n \frac1n
n^{p/{p^{'}}}\Big{(}\int_{\mathbb{H}^{'}}\Big{\|}\int_{\mathbb{H}}T\mathbf{f}(s_i+h)\sigma(h)\,
d_{\mu_{\mathbb{H}}}(h)\Big{\|}^{p^{'}}d_{{\mu}_{\mathbb{H}^{'}}}(\sigma+\mathbb{H}^{\bot})\Big{)}^{p/{p^{'}}}\Big{]}^{1/p}\\
&\leq \Big{[}\sum_{i=1}^n \frac1n
n^{p/{p^{'}}}\|\,T\,|\mathcal{FT}_p^{\mathbb{H}}\|^p
\Big{(}\int_{\mathbb{H}}
\|\mathbf{f}(s_i+h)\|^{p}\,d_{\mu_{\mathbb{H}}}(h)\Big{)}\Big{]}^{1/p}\\
&= \,n^{1/{p^{'}}}
\|\,T\,|\mathcal{FT}_p^{\mathbb{H}}\|\|\mathbf{f}|L_p(\mathbb{G})\|\mbox{.}
\\
\end{split}
\end{equation*}
We used Minkowski's inequality in the third line.
\end{proof}

\begin{prop}
\label{CGZ} If $\mathbb{F}$ is an infinite, compact and connected
LCA group then for any nonnegative integers $a\,,\,b$, there exist
$c(a,b)\,,\,C(a,b)
> 0 $ such that $$c(a,b)\,\|\,T\,|\mathcal{FT}_p^{\mathbb{T}}\|\leq
\|\,T\,|\mathcal{FT}_p^{\mathbb{R}^a\times\mathbb{Z}^b\times\mathbb{F}}\|\leq
C(a,b)\,\|\,T\,|\mathcal{FT}_p^{\mathbb{T}}\|.$$ Therefore we have
$\mathcal{FT}_p^{\mathbb{R}^a\times\mathbb{Z}^b\times\mathbb{F}}
\,=\,\mathcal{FT}_p $.
\end{prop}
\begin{proof}
If $a\,=\,b\,=\,0$ then by applying ($\ref{FZ}$) \,we have
$\|\,T'\,|\mathcal{FT}_p^{\mathbb{F}^{'}}\|\,=\,\|\,T'\,|\mathcal{FT}_p^{\mathbb{Z}}\|$
where $\mathbb{F}^{'}$ is the dual group of $\mathbb{F}$ and hence
a nontrivial torsion-free group with the discrete topology.
Therefore it follows that
\begin{equation}
\begin{split}
\label{CCG}
\|\,T\,|\mathcal{FT}_p^{\mathbb{F}}\|\,&=\,\|\,T'\,|\mathcal{FT}_p^{\mathbb{Z}}\|\,=\,\|\,T\,|\mathcal{FT}_p^{\mathbb{T}}\| \quad \mbox{ and} \\
 \mathcal{FT}_p^{\mathbb{F}}&=\mathcal{FT}_p. \\
\end{split}
\end{equation}
If $a\,+\,b\,\geq\,1$ then by applying Theorem $\ref{HNG}$\ and
Corollary $\ref{NRTZG}$, we have
$$\mathcal{FT}_p^{\mathbb{R}^a\times\mathbb{Z}^b\times\mathbb{F}}=\mathcal{FT}_p^{\mathbb{T}\times\mathbb{F}}$$
and in fact
$$c(a,b)\,\|\,T\,|\mathcal{FT}_p^{\mathbb{T}\times\mathbb{F}}\|\leq
\|\,T\,|\mathcal{FT}_p^{\mathbb{R}^a\times\mathbb{Z}^b\times\mathbb{F}}\|\leq
C(a,b)\,\|\,T\,|\mathcal{FT}_p^{\mathbb{T}\times\mathbb{F}}\|$$
for some positive reals $c(a,b)\,,\,C(a,b)$. Here
$c(a,b)=(\frac2{\pi})^{ab}$, $C(a,b)=1$. Note that
$\mathbb{T}\times\mathbb{F}$ is connected and compact. Thus we
again have
$$\|\,T\,|\mathcal{FT}_p^{\mathbb{T}\times\mathbb{F}}\|=\|\,T\,|\mathcal{FT}_p^{\mathbb{T}}\|
\mbox{ \ by (\ref{CCG})}.$$
\end{proof}

If $\mathbb{F}$ has only finitely many components and $\mathbb{C}$
is the component(maximal connected set) of the identity element
then $\mathbb{F}/ \mathbb{C}$ is finite and we have the following
result.

\begin{thm}
\label{FIN} Let $\mathbb{F}$ be an infinite compact LCA group with
$n$ components. Then for any nonnegative integers $a\,,\,b$, \
there exists $c(a,b)\,,\,C(a,b)
> 0 $ such that $$c(a,b)\,\|\,T\,|\mathcal{FT}_p^{\mathbb{T}}\|\leq
\|\,T\,|\mathcal{FT}_p^{\mathbb{R}^a\times\mathbb{Z}^b\times\mathbb{F}}\|\leq
n^{1/p^{'}}C(a,b)\,\|\,T\,|\mathcal{FT}_p^{\mathbb{T}}\|.$$ And
therefore
$\mathcal{FT}_p^{\mathbb{R}^a\times\mathbb{Z}^b\times\mathbb{F}}
\,=\,\mathcal{FT}_p $.
\end{thm}
\begin{proof}
Let $\mathbb{C}$ be the component of the identity element of
$\mathbb{F}$ then the factor group
$(\mathbb{R}^a\times\mathbb{Z}^b\times\mathbb{F})/(\mathbb{R}^a\times\mathbb{Z}^b\times\mathbb{C})$
is a finite group with $n$ elements. Thus by Lemma \ref{FCC},
$$\|\,T\,|\mathcal{FT}_p^{\mathbb{R}^a\times\mathbb{Z}^b\times\mathbb{C}}\|\leq
\|\,T\,|\mathcal{FT}_p^{\mathbb{R}^a\times\mathbb{Z}^b\times\mathbb{F}}\|
\leq
n^{1/p^{'}}\|\,T\,|\mathcal{FT}_p^{\mathbb{R}^a\times\mathbb{Z}^b\times\mathbb{C}}\|
$$ and
$$\mathcal{FT}_p^{\mathbb{R}^a\times\mathbb{Z}^b\times\mathbb{F}}=
\mathcal{FT}_p^{\mathbb{R}^a\times\mathbb{Z}^b\times\mathbb{C}}.$$
Then the statement follows from Proposition \ref{CGZ}.
\end{proof}

This Theorem can not say anything about the case when $\mathbb{F}$
has infinitely many components, because of the factor
$n^{1/p^{'}}$.

Theorem $\ref{FIN}$ can be extended when our scope goes beyond the
boundary of compactly generated LCA groups.
\begin{thm}
\label{GFE} Let $\mathbb{E}$ be a nontrivial torsion free group
with the discrete topology and $\mathbb{F}$ an infinite, compact
and connected LCA group. Then for any LCA group $\mathbb{G}$,\
$$\|\,T\,|\mathcal{FT}_p^{\mathbb{E}\times\mathbb{G}}\|=\|\,T\,|\mathcal{FT}_p^{\mathbb{Z}\times\mathbb{G}}\|$$
and
$$\|\,T\,|\mathcal{FT}_p^{\mathbb{F}\times\mathbb{G}}\|=\|\,T\,|\mathcal{FT}_p^{\mathbb{T}\times\mathbb{G}}\|.$$

If $\tilde{\mathbb{F}}$ is a compact LCA group with $n$ components
then for any nonnegative integer $k$ there exist \ $c(k), C(k)$
such that \
$$c(k)\,\|\,T\,|\mathcal{FT}_p^{\mathbb{T}\times\mathbb{G}}\|\leq
\|\,T\,|\mathcal{FT}_p^{\mathbb{R}^k\times\mathbb{E}^l\times\tilde{\mathbb{F}}^m\times\mathbb{G}}\|
\leq
n^{m/p^{'}}C(k)\,\|\,T\,|\mathcal{FT}_p^{\mathbb{T}\times\mathbb{G}}\|$$
and hence
$$\mathcal{FT}_p^{\mathbb{R}^k\times\mathbb{E}^l\times\tilde{\mathbb{F}}^m\times\mathbb{G}}=
\mathcal{FT}_p^{\mathbb{T}\times\mathbb{G}}$$ where \,
$l,\,m\,=\,0 \ \mbox{or} \ 1 \ \mbox{and} \ k+l+m \geq 1$. In
particular if $\mathbb{G}$ is the trivial group then
$$\mathcal{FT}_p^{\mathbb{R}^k\times\mathbb{E}^l\times\tilde{\mathbb{F}}^m}=\mathcal{FT}_p.$$
\end{thm}
\begin{proof}
First, $\mathbb{Z}$ is isomorphic with an open subgroup of
$\mathbb{E}$, and then $\mathbb{Z}\times\mathbb{G}$ is an open
subgroup of $\mathbb{E}\times\mathbb{G}$. By applying (\ref{OS})
we have
$\|\,T\,|\mathcal{FT}_p^{\mathbb{Z}\times\mathbb{G}}\|\leq\|\,T\,|\mathcal{FT}_p^{\mathbb{E}\times\mathbb{G}}\|$.
Conversely, for any simple $X$-valued function $\mathbf{f}$ which
is defined on $\mathbb{E}\times\mathbb{G}$ and has finite
$L_p$-norm, the support of $\mathbf{f}$ is a subset of
$\mathbb{Z}^k\times\mathbb{G}$ for some positive integer $k$.
Therefore by Lemma $\ref{ZNG}$ it follows that
$\|\,T\,|\mathcal{FT}_p^{\mathbb{E}\times\mathbb{G}}\|\leq\sup_k\|\,T\,|\mathcal{FT}_p^{\mathbb{Z}^k\times\mathbb{G}}\|
=\|\,T\,|\mathcal{FT}_p^{\mathbb{Z}\times\mathbb{G}}\|$.

Second, from the duality argument it follows that
$$\|\,T\,|\mathcal{FT}_p^{\mathbb{F}\times\mathbb{G}}\|=\|\,T^{'}\,|\mathcal{FT}_p^{\mathbb{F}^{'}\times\mathbb{G}^{'}}\|=
\|\,T^{'}\,|\mathcal{FT}_p^{\mathbb{Z}\times\mathbb{G}^{'}}\|=\|\,T\,|\mathcal{FT}_p^{\mathbb{T}\times\mathbb{G}}\|.$$

The rest follows by applying the above two results, Proposition
\ref{EQRZ}, \ref{EQRT} and  Lemma \ref{FCC}.
\end{proof}

\begin{rem}
If $\mathbb{G}$ is an LCA group and $\mathbb{C}$ is the component
of the identity of $\mathbb{G}$ then $\mathbb{C}$ is a closed
normal subgroup of  $\mathbb{G}$ and the factor group $\mathbb{G}/
\mathbb{C}$ is totally disconnected and Hausdorff. From (24.45) of
\cite{Hw}, if $\mathbb{C}$ is open then $\mathbb{G}$ is isomorphic
with $\mathbb{C}\times(\mathbb{G}/\mathbb{C})$. In particular
$\mathbb{C}$ is open when $\mathbb{G}$ is locally connected. Hence
to solve the question, denoted by (P), whether
$\|\,T\,|\mathcal{FT}_p^{\mathbb{G}}\|$ is equivalent to
$\|\,T\,|\mathcal{FT}_p^{\mathbb{T}}\|$ or not, it is useful to
find $\mathbb{C}$ and $\mathbb{G}/\mathbb{C}$. Since $\mathbb{C}$
is isomorphic with $\mathbb{R}^n\times \mathbb{K}$, where $n$ is a
nonnegative integer and $\mathbb{K}$ is a compact connected group,
see Theorem 9.14 of\cite{Hw}, the answer for (P) is affirmative
when $\mathbb{G}/\mathbb{C}$ is good.
\end{rem}

\begin{thm}
Let $\mathbb{G}$ be an LCA group with $n$ components. Then there
are positive real numbers $c$ and $C$ such that
\begin{equation}
\label{FF}
c\,\|\,T\,|\mathcal{FT}_p^{\mathbb{T}}\|\leq\|\,T\,|\mathcal{FT}_p^{\mathbb{G}}\|\leq
n^{1/p^{'}}C\,\|\,T\,|\mathcal{FT}_p^{\mathbb{T}}\|
\end{equation}
and hence $\mathcal{FT}_p^{\mathbb{G}}=\mathcal{FT}_p$.
\end{thm}
\begin{proof}
Let $\mathbb{C}$ be the component of the identity,
$\mathbb{G}/\mathbb{C}$ is a finite LCA group. Hence by Lemma
$\ref{FCC}$ we have
$\|\,T\,|\mathcal{FT}_p^{\mathbb{C}}\|\leq\|\,T\,|\mathcal{FT}_p^{\mathbb{G}}\|\leq
n^{1/p^{'}}\|\,T\,|\mathcal{FT}_p^{\mathbb{C}}\|$. And
$\mathbb{C}$ is isomorphic with $\mathbb{R}^k\times\mathbb{K}$,
where $k$ is a nonnegative integer and $\mathbb{K}$ is a compact
connected group. $\mathbb{K}$ is trivial or infinite. If
$\mathbb{K}$ is trivial then by applying Proposition \ref{EQRT} it
follows that
$$\,c\,\|\,T\,|\mathcal{FT}_p^{\mathbb{T}^k}\|\leq\|\,T\,|\mathcal{FT}_p^{\mathbb{C}}\|\leq
C\,\|\,T\,|\mathcal{FT}_p^{\mathbb{T}^k}\|$$ for some $c,\ C >0$.
And by  ($\ref{CCG}$),
$\|\,T\,|\mathcal{FT}_p^{\mathbb{T}^k}\|=\|\,T\,|\mathcal{FT}_p^{\mathbb{T}}\|$.
Thus we have
\begin{equation}
\label{FFF}
c\,\|\,T\,|\mathcal{FT}_p^{\mathbb{T}}\|\leq\|\,T\,|\mathcal{FT}_p^{\mathbb{C}}\|\leq
C\,\|\,T\,|\mathcal{FT}_p^{\mathbb{T}}\|. \end{equation} If
$\mathbb{K}$ is infinite, we also have (\ref{FFF}) by Proposition
\ref{CGZ}. Therefore the inequality ($\ref{FF}$)
 follows, and $\mathcal{FT}_p^{\mathbb{G}}=\mathcal{FT}_p$.
\end{proof}

\begin{thm}
\label{GFF} If $\mathbb{C}$ is the component of the identity of an
LCA group $\mathbb{G}$ and $\mathbb{G}/\mathbb{C}$ is a
(nontrivial torsion free group)$\times$(finite group with cardinal
$n$)\ with the discrete topology, then
$\|\,T\,|\mathcal{FT}_p^{\mathbb{C}}\|
\leq\|\,T\,|\mathcal{FT}_p^{\mathbb{G}}\| \leq
n^{1/{p^{'}}}\|\,T\,|\mathcal{FT}_p^{\mathbb{C}\times\mathbb{Z}}\|$
and $\mathcal{FT}_p^{\mathbb{G}}=\mathcal{FT}_p$.
\end{thm}
\begin{proof}
We have that $\mathbb{G}$ is isomorphic with $\mathbb{C}\times
(\mathbb{G}/\mathbb{C})$ because $\mathbb{C}$ is open. We can
apply Proposition \ref{OPP} and Theorem \ref{GFE} to obtain this
theorem.
\end{proof}

As for the class of Banach spaces which has Fourier type $p$ for
an infinite LCA group, we have the following:
\begin{cor}
\rm{i)} Under the same conditions as Theorem \ref{GFE} we have
$$\mathrm{FT}_p^{\mathbb{R}^k\times\mathbb{E}^l\times\tilde{\mathbb{F}}^m}=\mathrm{FT}_p.$$
\rm{ii)} Under the same conditions as Theorem \ref{GFF},
$$\mathrm{FT}_p^{\mathbb{G}}=\mathrm{FT}_p.$$
\end{cor}

\textbf{Acknowledgements.} The author would like to thank
Professor Chansun Choi and a colleague Hun Hee Lee for comments.

\bibliographystyle{amsplain}

\begin{thebibliography}{1}

\bibitem{An}
M. E. Andersson, \textit{On the vector valued Hausdorff-Young
inequality}, Ark. Mat. \textbf{36}(1998), no.1, 1-30.

\bibitem{Br}
J. Bourgain, \textit{Vector-valued Hausdorff-Young inequalities
and applications, in Geometric Aspects of Fucntional Analysis},
Lecture Notes in Math. 1317, p239-249, Springer-Verlag,
Berlin-Heidelberg, 1988.

\bibitem{Fol}
G. B. Folland, \textit{A course in abstract harmonic analysis},
(Studies in Advanced Mathematics) CRC Press, Boca Raton, FL, 1995.

\bibitem{Gar}
J. Garcia-Cuerva; K. S. Kazaryan; V. I. Kolyada; J. L. Torrea,
\textit{The Hausdorff-Young inequality with vector-valued
coefficients and applications}, translation in Russian Math.
Surveys \textbf{53} (1998), no. 3, 435-513.

\bibitem{Hw}
E. Hewitt; K. A. Ross, \textit{Abstract harmonic analysis}, Vol.
I. second edition. Grundlehren der Mathematischen Wissenschaften
[Fundamental Principles of Mathematical Sciences], 115.
Springer-Verlag, Berlin-New York, 1979.

\bibitem{Jo}
Frank Jones, \textit{Lebesque Integration on Euclidean Space},
Jones and Bartlett, 1993.

\bibitem{Ko}
H. K\"{o}nig, \textit{On the Fourier-coefficients of vector-valued
functions}, Math. Nachr. \textbf{152}(1991), 215-227.

\bibitem{Kw}
S. Kwapi\'{e}n, \textit{Isomorphic characterizarions of inner
product spaces by orthogonal series with vector-valued
coefficients}, Stud.Math. \textbf{44} (1972) 583-595

\bibitem{Pee}
J. Peetre, \textit{Sur la transformation de Fourier des fonctions
a valeurs vectorielles}, (French) Rend. Sem. Mat. Univ. Padova
\textbf{42}(1969), 15-26.

\bibitem{Pie}
A. Pietsch; J. Wenzel, \textit{Orthonormal systems and Banach
space geometry}, Encyclopedia of Mathematics and its Applications
\textbf{70}. Cambridge University Press, Cambridge, 1998.

\bibitem{Re}
H. Reiter, \textit{Classical Harmonic Analysis and Locally compact
Groups}, Oxford U. Press, Oxford, 1968.

\bibitem{Rud}
W. Rudin, \textit{Fourier analysis on groups}, Interscience Tracts
in Pure and Applied Mathematics, No. 12 Interscience Publishers (a
division of John Wiley and Sons), New York-London 1962.

\end{thebibliography}
\providecommand{\bysame}{\leavevmode\hbox
to3em{\hrulefill}\thinspace}

\end{document}